\documentclass[a4paper,12pt]{amsart}
\usepackage{ifthen}
\usepackage{graphicx}
\usepackage{epstopdf}
\usepackage{mathrsfs}
\usepackage{color}
\nonstopmode \numberwithin{equation}{section}
\newtheorem{thm}{Theorem}
\newtheorem{cor}{Corollary}
\newtheorem{lem}{Lemma}

\theoremstyle{definition}
\newtheorem{defn}{Definition}

\newtheorem{rem}{Remark}

\newenvironment{pf}[1][]{%
 \vskip 3mm
 \noindent
 \ifthenelse{\equal{#1}{}}%
  {{\slshape Proof. }}%
  {{\slshape #1.} }%
 }%
{\qed\bigskip}

\newcounter{alphabet}
\newcounter{tmp}

\newcommand{\A}{{\mathcal A}}

\newcommand{\C}{{\mathbb C}}
\newcommand{\D}{{\mathbb D}}

\newcommand{\K}{{\mathcal K}}

\newcommand{\Q}{{\mathcal Q}}
\newcommand{\R}{{\mathbb R}}

\newcommand{\es}{{\mathcal S}}

\newcommand{\sphere}{{\widehat{\mathbb C}}}

\renewcommand{\Im}{{\,\operatorname{Im}\,}}

\renewcommand{\Re}{{\,\operatorname{Re}\,}}

\newcommand{\inv}{^{-1}}

\newcommand{\Gauss}{{\null_2F_1}}

\renewcommand{\arg}{\,{\operatorname{arg}\,}}

\newcommand{\aand}{{\quad\text{and}\quad}}
\newcommand{\sgn}{{\operatorname{sgn}\,}}

\begin{document}
\bibliographystyle{amsplain}
\title{Geometric properties of the shifted hypergeometric functions}

%\begin{center}
%{\tiny \texttt{FILE:~\jobname .tex,
%        printed: \number\year-\number\month-\number\day,
%        \thehours.\ifnum\theminutes<10{0}\fi\theminutes}
%}
%\end{center}
\author[T. Sugawa]{Toshiyuki Sugawa}
\address{Graduate School of Information Sciences, \\
Tohoku University, \\
Aoba-ku, Sendai 980-8579, Japan} \email{sugawa@math.is.tohoku.ac.jp}
\author[L.-M.~Wang]{Li-Mei Wang}
\address{School of Statistics,
University of International Business and Economics, No.~10, Huixin
Dongjie, Chaoyang District, Beijing 100029, China}
\email{wangmabel@163.com}

\keywords{subordination, hypergeometric functions,
strongly starlike function, spirallike function}
\subjclass[2010]{Primary 30C45; Secondary 33C05}
\begin{abstract}
We will provide sufficient conditions for the shifted hypergeometric function $z_2F_1(a,b;c;z)$
to be a member of a specific subclass of starlike functions in terms of the complex parameters
$a,b$ and $c.$
For example, we study starlikeness of order $\alpha,$ $\lambda$-spirallikeness of order $\alpha$
and strong starlikeness of order $\alpha.$
In particular, those properties lead to univalence of the shifted hypergeometric functions on the unit disk.
\end{abstract}
\maketitle

\section{Introduction}
The (Gaussian) hypergeometric functions appear in various
areas in Mathematics, Physics and Engineering and have proved to be
quite useful in many respects. Many of the common mathematical functions 
can be expressed in terms of  hypergeometric functions, or suitable limits of them. 
Their geometric properties in complex domains were, however, studied
only recently (in comparison with its long history). 
For instance, starlikeness and convexity are investigated in 1960's and afterwards.
See \cite[\S 4.5]{MM:ds}, K\"ustner \cite{Kus02}, \cite{Kus07},
H\"asto, Ponnusamy and Vuorinen \cite{HPV10} and the references therein. 
This sort of research is not only interesting in the viewpoint of classical analysis,
but also applicable in the theory of function spaces, integral transforms, 
convolutions and so on (see \cite{BNP06}, \cite{BPV02} and \cite{HPV10} for example).
It should also be noted that most of known results are restricted to real parameter cases.
In \cite{SW16sp}, the authors gave several sufficient conditions
for spirallikeness and strong starlikeness of the shifted
hypergeometric functions with complex parameters.
In the present note, we will extend the results for
general classes of starlike functions.

Let $\A$ denote the set of analytic functions on the unit disk
$\D=\{z\in\C: |z|<1\}$ in the complex plane
and consider the subclasses $\A_0=\{\varphi\in\A: \varphi(0)=1\}$
and $\A_1=\{f\in\A: f(0)=f'(0)-1=0\}$ of $\A.$
%We denote by $\es$ the subset of $\A_1$ consisting of
%univalent functions on $\D.$
A function $f\in\A_1$ is called {\it starlike} if $f$ is univalent and the image
$f(\D)$ is stalike with respect to the origin.
This property is characterized by the condition $\Re[zf'(z)/f(z)]>0$ on $\D.$
Robertson \cite{Rob36} refined this notion as follows.
For a constant $\alpha\in[0,1)$, a function $f\in\A_1$ is
called {\it starlike of order $\alpha$} if
$$
\Re\frac{zf'(z)}{f(z)}>\alpha,\quad z\in\D.
$$
As another refinement, Stankiewicz \cite{Stank66} and Brannan and Kirwan \cite{BK69}
introduced independently the class of strongly starlike functions of order $\alpha$ for $0<\alpha<1,$
which is defined by the condition $|\arg[zf'(z)/f(z)]|<\pi\alpha/2$ on $\D.$
For several geometric characterizations of this class, see \cite{SugawaDual}.
It is also known that a strongly starlike function has a quasiconformal extension to the complex plane
(see \cite{FKZ76}).

As an extension of starlikeness, spirallikeness is natural and useful.
For a real number $\lambda$ with $|\lambda|<\pi/2,$
a function $f\in\A_1$ is called {$\lambda$-spirallike} if $f$ is univalent and if
the $\lambda$-spiral of the form $w\exp(-te^{i\lambda}),~0\le t<+\infty,$ is contained
in $f(\D)$ for each $w\in f(\D).$
This is characterized by the condition $\Re[e^{-i\lambda}zf'(z)/f(z)]>0$ on $\D.$
Libera \cite{Libera67} refined this notion as follows:
a function $f\in\A_1$ is called {\it $\lambda$-spirallike of order $\alpha$} if
$$
\Re\left(e^{-i\lambda}\frac{zf'(z)}{f(z)}\right)>\alpha\cos\lambda,\quad z\in\D.
$$
%Let $\spl(\alpha,\lambda)$ denote the class of $\lambda$-spirallike functions of order $\alpha$.

Ma and Minda \cite{MM92B} proposed a unifying way to treat these classes as follows.
Let
$$
\es^*(\varphi)=\{f\in\A_1: zf'/f\prec\varphi\}.
$$
Here, $\varphi\in\A_0$ and
the symbol $g\prec h$ means subordination,
that is, $g=h\circ\omega$ for an analytic function $\omega$
on $\D$ with $|\omega(z)|\le|z|$ for $z\in\D.$
If $h$ is univalent in $\D$, $g\prec h$
if and only if $g(0)=h(0)$ and $g(\D)\subset h(\D)$.
We summarize typical choices of $\varphi$ and the corresponding classes
$\es^*(\varphi)$ in Table 1.
See, for instance, \cite{Goodman:univ} for details about the special classes
of univalent functions.

\renewcommand\arraystretch{3}
\begin{table}[h]
\begin{tabular}{|c|c|} \hline
$\varphi(z)$ & $\es^*(\varphi)$ \\ \hline\hline
$\varphi_0(z)=\dfrac{1+z}{1-z}$ & starlike functions \\ \hline
$\varphi_\alpha(z)=\dfrac{1+(1-2\alpha)z}{1-z}$ &  starlike functions of order $\alpha$ \\ \hline
$\phi_\alpha(z)=\left(\dfrac{1+z}{1-z}\right)^\alpha$ & strongly starlike functions of order $\alpha$ \\ \hline
$\psi_\lambda(z)=\dfrac{1+e^{2i\lambda}z}{1-z}$ & $\lambda$-spirallike functions \\ \hline
$\psi_{\lambda,\alpha}(z)=\dfrac{1+[e^{2i\lambda}-\alpha(1+e^{2i\lambda})]z}{1-z}$
& $\lambda$-spirallike functions of order $\alpha$\\ \hline
\end{tabular}
\medskip
\caption{$\varphi$ and the corresponding class $\es^*(\varphi)$}
\end{table}

The main aim in the present paper is to give a sufficient condition
for the shifted hypergeometric function $f(z)=z\Gauss(a,b;c;z)$ to be a member
of the class $\es^*(\varphi)$ for a certain $\varphi.$
To state our main theorem, we introduce a class of admissible functions $Q=\varphi-1.$

\begin{defn}\label{def:Q}
We denote by $\Q$ the set of all analytic functions $Q$ satisfying the following.
There exist finitely many points $\zeta_1,\dots,\zeta_n$ in $\partial\D$ with
the following five conditions:
\renewcommand{\labelenumi}{({\roman{enumi}})}
\begin{enumerate}
\item
$Q$ is analytic and univalent on a neighbourhood of
$\overline\D\setminus\{\zeta_1,\dots,\zeta_n\}.$
\item
$Q(0)=0.$
\item
the function $Q(z)$ has a limit $w_j$ in $\sphere$
as $z\to\zeta_j$ in $\D$ for each $j=1,\dots,n.$
\item\label{item:fin}
When $w_j\in\C,$  $(Q(z)-w_j)^{\beta_j}$ extends to a univalent function
on a neighbourhood of $z=\zeta_j$ for some number $\beta_j>1.$
\item
When $w_j=\infty,$  $Q(z)^{-\beta_j}$ extends to a univalent function
on a neighbourhood of $z=\zeta_j$ for some number $\beta_j\ge 1.$
\item\label{item:P}
When $w_j=\infty$ and $\beta_j=1,$ the derivative of $P=1/Q$ satisfies
$\zeta_jP'(\zeta_j)\in\C\setminus[0,1].$
\end{enumerate}
\end{defn}

We note that a similar class was introduced by Miller and Mocanu
(see Definition 2.2b in \cite{MM:ds}) but our class is more restrictive.
It is easy to check the above conditions for all the functions in Table 1
with $0\le \alpha<1$ and $-\pi/2<\lambda<\pi/2.$
As an example, we examine the function $Q=\psi_{\lambda,\alpha}-1.$
In this case $n=1$ and $\zeta_1=1, w_1=\infty, \beta_1=1.$
It is enough to check condition (vi) because the other ones are, more or less, obvious.
Let $P(z)=1/Q(z)=(1-z)/[(1-\alpha)(1+e^{2i\lambda})z].$
Then $\zeta_1P'(\zeta_1)=P'(1)=-1/[(1-\alpha)(1+e^{2i\lambda})].$
If $P'(1)\in[0,1],$ then $\lambda=0$ and $-1/(1-\alpha)\in[0,2],$ which is impossible
for $0\le \alpha<1.$
Thus, we have seen that $Q\in\Q$ in this case.

The following is our main theorem, from which we will derive several consequences
in Sections 3 and 4.

\begin{thm}\label{thm:main}
Let $a,b,c\in\C$ with $c\ne0,-1,-2,\dots$ and let $\varphi=1+Q$ for some $Q\in\Q.$
Then $f(z)=z\Gauss(a,b;c;z)$ belongs to $\es^*(\varphi)$ if
\begin{equation}\label{eq:main1}
0<-2\Re\big[(pQ(\zeta)+ab)\overline{\zeta Q'(\zeta)}\big]
\quad\text{and}
\end{equation}
\begin{equation}\label{eq:main2}
|B(\zeta)|^2-|A(\zeta)|^2\leq -2\Re\big[(pQ(\zeta)+ab)\overline{\zeta Q'(\zeta)}\big]
\end{equation}
for $\zeta\in \partial\D\setminus\{\zeta_1,\dots,\zeta_n\},$
where $p=a+b+1-c,~ A(z)=Q(z)(Q(z)+c-1)$ and
$B(z)=(Q(z)+a)(Q(z)+b).$
\end{thm}

In practical computations, it is convenient to express $|B|^2-|A|^2$ in the following
form:
\begin{align}\label{eq:AB}
|B(z)|^2-|A(z)|^2
&=|w|^2(2\Re[p\bar w]+|a|^2+|b|^2-|c-1|^2) \\
&\quad +(2\Re[a\bar w]+|a|^2)(2\Re[b\bar w]+|b|^2),
\notag
\end{align}
where $w=Q(z).$

\begin{rem}\label{rem:AB}
The condition \eqref{eq:main1} can be weakened to
$0\le -2\Re[(pQ(\zeta)+ab)\overline{\zeta Q'(\zeta)}]$ if, instead, the condition
$A(\zeta)\ne B(\zeta)$ is guaranteed.
\end{rem}

\begin{rem}\label{rem:convex}
We can also obtain a convexity counterpart as follows.
Let $\K(\varphi)$ be the class of functions $f\in\A_1$
satisfying $1+zf''/f'\prec\varphi$ for a given $\varphi\in \A_0.$
For $f\in\A_1,$ as is well known, $f\in\K(\varphi)$ if and only if
$zf'\in\es^*(\varphi).$
When $f(z)=\frac{c}{ab}(\Gauss(a,b,c;z)-1),$ by \eqref{eq:der} below,
we have $zf'(z)=z\Gauss(a+1,b+1,c+1;z).$
Therefore, we have a sufficient condition for the function
$f(z)=\frac{c}{ab}(\Gauss(a,b,c;z)-1)$ to be a member of $\K(\varphi)$
as an immediate consequence of Theorem \ref{thm:main}, though we do not
state it separately.
\end{rem}

\section{Preliminaries and proof of the main theorem}

First we recall a definition and basic properties of hypergeometric functions.
The (Gaussian) {\it hypergeometric function} $\Gauss(a,b;c;z)$
with complex parameters $a,b,c$~ $(c\not=0,-1,-2,\dots)$
is defined by the power series
$$
\Gauss(a,b;c;z)=\sum_{n=0}^{\infty}\frac{(a)_n(b)_n}{(c)_nn!}z^n
$$
on $|z|<1,$ where $(a)_n$ is the {\it Pochhammer symbol}; namely, $(a)_0=1$ and
$(a)_n=a(a+1)\cdots(a+n-1)=\Gamma(a+n)/\Gamma(a)$ for $n=1,2,\dots.$
We note that $\Gauss(a,b;c;z)$ analytically continues on
the slit plane $\C\setminus[1,+\infty).$
Note here that the hypergeometric function is symmetric in the parameters $a$ and $b$
in the sense that $\Gauss(a,b;c;z)=\Gauss(b,a;c;z).$
As is well known, the hypergeometric function $F(z)=\Gauss(a,b;c;z)$
is characterized as the solution to the hypergeometric differential equation
\begin{equation}\label{eq:Gauss}
(1-z)zF''(z)+[c-(a+b+1)z]F'(z)-abF(z)=0
\end{equation}
with the initial condition $F(0)=1.$
We also note the derivative formula:
\begin{equation}\label{eq:der}
\frac{d}{dz}\Gauss(a,b;c;z)=\frac{ab}{c}\Gauss(a+1,b+1;c+1;z).
\end{equation}
For more properties of hypergeometric functions, we refer to
\cite{AS:hand}, \cite{Temme:sp} and \cite{WW:anal} for example.

As in \cite{Lewis79}, \cite{SW16sp} or in the proof of \cite[Theorem 2.12]{Rus:conv}),
our proof of the main theorem will be based on the following (easier) variant of Julia-Wolff theorem
(see \cite[Prop.~4.13]{Pom:bound} for instance) and
the hypergeometric differential equation \eqref{eq:Gauss}.

\begin{lem}\label{lem:Jack}
Let $z_0\in\C$ with $|z_0|=r\ne0$ and
let $\omega$ be a non-constant analytic function on a neighbourhood of
$\{z: |z|<r\}\cup\{z_0\}$ with $\omega(0)=0$.
If $|\omega(z)|\le|\omega(z_0)|$ for $|z|<r,$
then $z_0\omega'(z_0)=k\omega(z_0)$ for some $k\geq 1$.
\end{lem}

We are now ready to prove our main theorem.

\begin{pf}[Proof of Theorem \ref{thm:main}]
Let $\zeta_j, w_j$ and $\beta_j,~j=1,\dots,n,$ be as in Definition \ref{def:Q} and
set $F(z)=\Gauss(a,b;c;z)$ and $f(z)=zF(z).$
Let us try to show that $zf'(z)/f(z)=1+zF'(z)/F(z)\prec \varphi(z)=1+Q(z).$
Put $q(z)=zF'(z)/F(z).$
Let $0<r\le 1$ be the largest possible number such that $q(z)\in \Omega:=Q(\D)$
for $|z|<r.$
We set $\omega(z)=Q\inv(q(z))$ for $|z|<r.$
Then $\omega(0)=0,$ $|\omega(z)|<1$ and $q(z)=Q(\omega(z))$ on $|z|<r.$
It thus suffices to show $r=1.$
Suppose, to the contrary, that $r<1.$
Then, there is a $z_0\in\C$ with $|z_0|=r$ such that $w_0:=q(z_0)\in\partial\Omega,$
where the boundary is taken in the Riemann sphere $\sphere=\C\cup\{\infty\}.$
We first assume that $w_0\ne w_j=Q(\zeta_j)$ for $j=1,\dots,n.$
Since $Q(z)$ is univalent near the point $z=\zeta_0:=Q\inv(w_0),$ the function $\omega(z)$
extends to $z=z_0$ analytically and satisfies $\omega(z_0)=\zeta_0\in\partial\D.$
Differentiating both sides of the relation $zF'(z)=F(z)Q(\omega(z)),$ we get
$$
zF''(z)=Q'(\omega(z))\omega'(z)F(z)+[Q(\omega(z))-1]F'(z).
$$
Substituting it into \eqref{eq:Gauss} and multiplying with $z,$ we obtain
$$
\big[c-(a+b+1)z+(1-z)(Q(\omega(z))-1)\big]zF'(z)=\big[(z-1)Q'(\omega(z))\omega'(z)+ab\big]zF(z).
$$
We replace $zF'(z)$ by $F(z)Q(\omega(z))$ in the last formula and rearrange it
to obtain
$$
\big[c-(a+b+1)z+(1-z)(Q(\omega(z))-1)\big]F(z)Q(\omega(z))=\big[(z-1)Q'(\omega(z))\omega'(z)+ab\big]zF(z).
$$
Since $F$ is not identically zero, we have
$$
\big[c-(a+b+1)z+(1-z)(Q(\omega(z))-1)\big]Q(\omega(z))=\big[(z-1)Q'(\omega(z))\omega'(z)+ab\big]z,
$$
which further leads to
\begin{equation}\label{eq:main}
Q(\omega)(Q(\omega)+c-1)+z\omega'Q'(\omega)
=\big[(Q(\omega)+a)(Q(\omega)+b)+z\omega'Q'(\omega)\big]z.
\end{equation}
Lemma \ref{lem:Jack} now implies that
$z_0\omega'(z_0)=k\omega(z_0)=k\zeta_0$ for some $k\geq 1.$
Letting $z=z_0$ in \eqref{eq:main} and using this result, we obtain
$$
w_0(w_0+c-1)+k\zeta_0Q'(\zeta_0)=\big[(w_0+a)(w_0+b)+k\zeta_0Q'(\zeta_0)\big]z_0.
$$
Let $A=w_0(w_0+c-1)$ and $B=(w_0+a)(w_0+b).$
In order to get a contradiction, it is enough to show the two
conditions:
\begin{enumerate}
\item[(I)]
$|A+k\zeta_0Q'(\zeta_0)|\ge|B+k\zeta_0Q'(\zeta_0)|,$
\item[(II)]
the two equalities $A+k\zeta_0Q'(\zeta_0)=0$ and $B+k\zeta_0Q'(\zeta_0)=0$ do not hold simultaneously.
\end{enumerate}
Since $B-A=(a+b+1-c)w_0+ab=pw_0+ab,$
condition (II) follows from the assumption \eqref{eq:main1}
(or instead the condition $A(\zeta)\ne B(\zeta)$ as is stated
in Remark \ref{rem:AB}).
Condition (I) means that $k\zeta_0Q'(\zeta_0)$ belongs
the half-plane $|w+A|\ge|w+B|.$
Note that the inequality $|w+A|^2\ge|w+B|^2$
is equivalent to $|A|^2-|B|^2\ge 2\Re[(B-A)\bar w].$
The assumptions \eqref{eq:main1} and \eqref{eq:main2}
imply now that $|A|^2-|B|^2\ge
2\Re[(B-A)\overline{\zeta_0Q'(\zeta_0)}]
\ge 2\Re[(B-A)\overline{k\zeta_0Q'(\zeta_0)}].$
Hence, condition (I) follows.
In this way, we have excluded
the possibility that $w_0\in\partial\Omega\setminus\{w_1,\dots,w_n\}.$

The remaining possibility is that $w_0=w_j$ for some $j.$
We first consider the case when $w_0=w_j\in\C.$
By a local property of analytic functions (see \cite[Chap.4, \S 3.3]{Ahlfors:ca}),
$q(z)=h(z)^m+w_j$ near $z=z_0,$ where $m$ is a positive integer and
$h(z)$ is a univalent analytic function near $z=z_0$ with $h(z_0)=0.$
In particular, the image of the disk $|z|<r$ under $q$ covers a (truncated) sector of opening
angle $\pi-\varepsilon$ with vertex at $w_j$ for an arbitrarily small $\varepsilon>0.$
On the other hand, condition (iv) implies that the interior angle of the domain
$\Omega=Q(\D)$ at $w_j$ is $\pi/\beta_j<\pi.$
Therefore, this case does not occur.
Next we consider the case when $\zeta_0=\zeta_j$ and $w_j=\infty.$
If $\beta_j>1,$ then the same argument as in the previous case works to
conclude that this is impossible.
Thus, $\beta_j=1.$
In this case, $P(z)=1/Q(z)$ is conformal at $z=\zeta_0,$ and therefore $P(\zeta_0)=0$
and $P'(\zeta_0)\ne0.$
Since $Q=1/P$ and $Q'=-P'/P^2,$ the formula \eqref{eq:main} turns to
$$
1+(c-1)P(\omega)-z\omega'P'(\omega)
=\big[(1+aP(\omega))(1+bP(\omega))-z\omega'P'(\omega)\big]z.
$$
We now let $z\to z_0$ to obtain further
$$
1-z_0\omega'(z_0)P'(\zeta_0)=\big[1-z_0\omega'(z_0)P'(\zeta_0)\big]z_0.
$$
In view of $|z_0|<1$ and $z_0\omega'(z_0)=k\omega(z_0)=k\zeta_0,$ we conclude that
$$
1-z_0\omega'(z_0)P'(\zeta_0)=1-k\zeta_0P'(\zeta_0)=0,
$$
which violates condition (vi).
Now all the possibilities have been excluded.
We thus conclude that $r=1$ as required.
\end{pf}

\section{Starlikeness and spirallikeness}

Note that $\varphi_\alpha=\psi_{0,\alpha}$ and $\psi_\lambda=\psi_{\lambda,0}$ in Table 1.
Thus the family $\psi_{\lambda,\alpha}$ covers the cases of  starlike functions of order $\alpha$
and $\lambda$-spirallike functions.
In order to apply Theorem \ref{thm:main} to the function $\psi_{\lambda,\alpha},$
we consider the function
$$
Q(z)=\psi_{\lambda,\alpha}(z)-1
=\dfrac{(1-\alpha)(1+e^{2i\lambda})z}{1-z}
=\dfrac{2\mu z}{1-z}
$$
for $\alpha\in[0,1)$ and $|\lambda|<\pi/2,$
where
$$
\mu=(1-\alpha)e^{i\lambda}\cos\lambda.
$$
Let $\zeta=e^{i\theta}\in \partial{\D}\setminus\{1\}$ and
$s=\cot \frac{\theta}{2}$ for $0<\theta<2\pi$.
Simple computations show that
\begin{align*}
%\zeta&=&e^{i\theta}=\frac{s^2-1}{s^2+1}+i\frac{2s}{s^2+1},\\
Q(\zeta)&=\mu\left(\frac{1+\zeta}{1-\zeta}-1\right)=\mu(-1+is)\quad  \text{and}\\
\zeta Q'(\zeta)&=\frac{2\mu \zeta}{(1-\zeta)^2}=\frac{2\mu}{(e^{i\theta/2}-e^{-i\theta/2})^2}
=-\frac{\mu(1+s^2)}{2}.
\end{align*}
%Let $w=Q(\zeta)=\mu(-1+is)$.
Now the condition \eqref{eq:main1} is equivalent to the inequality
$$
-2\Re\big[(pQ(\zeta)+ab)\overline{\zeta Q'(\zeta)}\big]
=(1+s^2)\Re\left[|\mu|^2(-1+is)p+ab\bar\mu\right]>0.
$$
Since $s$ can be any real number, this inequality forces that $p\in\R$
and \eqref{eq:main1} is further equivalent to $(1+s^2)\big(\Re[ab\bar\mu]- p|\mu|^2\big)>0.$
By \eqref{eq:AB}, we see that $|B(\zeta)|^2-|A(\zeta)|^2$ is a polynomial in $s$ and
$$
|B(\zeta)|^2-|A(\zeta)|^2
%=-2p(s\Im\mu+\Re\mu)|\mu|^2(s^2+1)
=-2p\Im\mu |\mu|^2s^3+O(s^2)\quad(s\to\pm\infty).
$$
Since the right-hand side of \eqref{eq:main2} is a polynomial in $s$ of degree 2,
we need the condition $p\Im\mu=0$ for the inequality \eqref{eq:main2} to hold
for all $s\in\R.$
Hence, Theorem \ref{thm:main}  works only when $\lambda=0$ or $p=0.$
In this case, the conditions \eqref{eq:main1} and \eqref{eq:main2} reduce to, respectively,
\begin{equation}\label{eq2:main1}
\Re[ab\bar\mu]-p|\mu|^2>0\quad \text{and}
\end{equation}
\begin{align}\label{eq2:main2}
&(1+s^2)\left(\Re[ab\bar\mu]-p|\mu|^2\right)
+(2p\Re\mu-|a|^2-|b|^2+|c-1|^2)|\mu|^2(1+s^2)\\
&-(2s\Im[a\bar\mu]-2\Re[a\bar\mu]+|a|^2)(2s\Im[b\bar\mu]-2\Re[b\bar\mu]+|b|^2)
\geq 0.
\notag
\end{align}

Letting $\lambda=0,$
we can prove the following theorem, which gives a sufficient condition
for the shifted hypergeometric function to be starlike of order $\alpha.$
Note that this is a natural generalization of \cite[Theorem 1.2]{SW16sp}.

\begin{thm}\label{thm:st}
Let $\alpha$ be a real constant with $0\leq \alpha<1$
and $a,b,c$ be complex numbers with $ab\ne0$ and $c\ne0,-1,-2,\dots.$
Then the shifted hypergeometric function $z\Gauss(a,b;c;z)$ is starlike
of order $\alpha$ if the following conditions are satisfied:
\renewcommand{\labelenumi}{({\roman{enumi}})}
\begin{enumerate}
\item
$p=a+b+1-c$ is a real number,
\item
$\Re[ab]>p(1-\alpha),$
\item
$L\ge0, N\ge0$ and $LN-M^2\ge0,$ where
\begin{align*}
L&=\frac{\Re[ab]}{1-\alpha}+p(1-2\alpha)-|a|^2-|b|^2+|c-1|^2-4\Im[a]\Im[b], \\
M%&=\Im[a](-2\Im[b]+|b|^2/(1-\alpha))+\Im[b](-2\Im[a]+|a|^2/(1-\alpha)) \\
&=\frac{\Im[ab(\bar a+\bar b-2+2\alpha)]}{1-\alpha}, \aand \\
N&=\frac{\Re[ab]}{1-\alpha}+p(1-2\alpha)-|a|^2-|b|^2+|c-1|^2 \\
&\qquad -\left(2\Re[a]-\frac{|a|^2}{1-\alpha}\right)\left(2\Re[b]-\frac{|b|^2}{1-\alpha}\right).
\end{align*}
\end{enumerate}
\end{thm}

There are several ways to express the coefficients $L,M$ and $N.$
To rephrase them, it is sometimes convenient to use the following elementary formulae:
\begin{equation}\label{eq:zw}
\Re[zw]=\Re[z\bar w]-2\Im z\Im w
=-\Re[z\bar w]+2\Re z\Re w\quad\text{for}~ z,w\in\C.
\end{equation}

\begin{pf}
For the choice $\lambda=0,$ we have $\mu=1-\alpha$ and $w=(1-\alpha)(-1+is)$ in the above observations.
For convenience, we write $a=a_1+ia_2$ and $b=b_1+ib_2.$
Substituting these, the left-hand side of \eqref{eq2:main2} can be expressed by
\begin{align*}
&\qquad \ (1-\alpha)^2\big\{\Re[ab]/(1-\alpha)-p+2p(1-\alpha)-|a|^2-|b|^2+|c-1|^2\big\}(1+s^2) \\
&\qquad -(1-\alpha)^2\big\{2a_2s-2a_1+|a|^2/(1-\alpha)\big\}\big\{2b_2s-2b_1+|b|^2/(1-\alpha)\big\} \\
&\quad =(1-\alpha)^2(Ls^2-2Ms+N).
\end{align*}
The above quadratic polynomial in $s$ is non-negative if and only if $L\ge0, N\ge0$ and $LN-M^2\ge0.$
Thus the assertion follows.
\end{pf}

Let $a=2(1-\alpha)$ in Theorem \ref{thm:st}.
Then $p=b-c+3-2\alpha$ should be real; in other words, $\Im b=\Im c.$
Thus, we can consider $f(z)=z\Gauss(2-2\alpha,b+is;c+is;z)$ for real $b,c,s.$
If, in addition, $c-b\ge 0$ and  $b+c>3,$
$\Re[a(b+is)]-p(1-\alpha)=(b+c-3+2\alpha)(1-\alpha)>0$ and
$L, M, N$ have the simple forms $L=N=2b+p(1-2\alpha)-4(1-\alpha)^2-b^2+(c-1)^2
=(c-b)(b+c-3+2\alpha), ~M=0$
so that $LN-M^2=L^2\ge0$ and $L=N\ge0.$
Therefore, as a consequence of Theorem \ref{thm:st}, we have the next result
due to Ruscheweyh \cite[Theorem 2.12, p. 60]{Rus:conv}.

\begin{cor}\label{cor:a=2}
Let $a,b,c,s$ be real numbers with $0<a\le 2$, $3\le b+c$ and $b\le c.$
Then the function $f(z)=z\Gauss(a,b+is;c+is;z)$ is starlike of order $1-a/2$.
In particular, $f$ is starlike.
\end{cor}

\begin{pf}
As is accounted above, the assertion follows from Theorem \ref{thm:st} when $3<b+c.$
When $b+c=3,$ we first apply the theorem to the function
$z\Gauss(a,b+\varepsilon+is;c+\varepsilon+is;z)$ for $\varepsilon>0$ and let $\varepsilon \to0.$
The assertion follows from the fact that the class $\es^{*}(\varphi_{\alpha})$ is compact
(see \cite{Pin68}).
\end{pf}

\begin{rem}\label{rem:ex}
The starlikeness of $f$ follows from \cite[Theorem 1.2]{SW16sp} when $a=2$
(see \cite[Corollary 4.1]{SW16sp}).
However, it seems that the above corollary does not follow from it for general $a\in(0,2).$
\end{rem}

We will next apply our main theorem to the case when $p=0$ and $-\pi/2<\lambda<\pi/2.$
Then we should eliminate $c$ by using the relation $c=a+b+1.$
Our goal is to show the following.

\begin{thm}\label{thm:spl}
Let $\lambda$ and $\alpha$ be real numbers with $|\lambda|<\pi/2$
and $\alpha\in[0,1).$
Let $a,b$ be complex numbers with $a+b\ne -1,-2,-3\dots.$
Then the shifted hypergeometric function $f(z)=z\Gauss(a,b;a+b+1;z)$ is
$\lambda$-spirallike of order $\alpha$ if the following two conditions are satisfied:
\renewcommand{\labelenumi}{({\roman{enumi}})}
\begin{enumerate}
\item
$\Re\big[e^{-i\lambda}ab\big]>0,$
\item
$L\ge0, N\ge0$ and $LN-M^2\ge0,$ where
\begin{align*}
L&=\Re\big[e^{-i\lambda}ab(2-\alpha+(1-\alpha)e^{-2i\lambda})\big], \\
M&=\Im\big[e^{-i\lambda}ab(\bar a+\bar b-(1-\alpha)(1+e^{-2i\lambda}))\big],
\aand \\
N&=\Re\big[e^{-i\lambda}ab(2\bar a+2\bar b+\alpha-(1-\alpha)e^{-2i\lambda})\big]
-|a|^2|b|^2/((1-\alpha)\cos\lambda).
\end{align*}
\end{enumerate}
\end{thm}

\begin{rem}
By \cite[Theorem 1.1 (iii)]{SW16sp}, the condition $\Re[e^{-i\lambda}ab]\ge0$
is necessary for $f$ to be $\lambda$-spirallike.
\end{rem}

\begin{pf}
For convenience, we put
$$
\nu=(1-\alpha)\cos\lambda.
$$
%so that $\mu=\nu e^{i\lambda}.$
Recalling $p=a+b+1-c=0$ and $w=\mu(-1+is),$
we see that \eqref{eq2:main1} is equivalent to
$$
\Re[ab(1-\alpha)e^{-i\lambda}\cos\lambda]
=\nu\Re[abe^{-i\lambda}]>0.
$$
Similarly, the left-hand side of \eqref{eq2:main2} is expressed by
\begin{align*}
&\quad \ (1+s^2)\big\{\Re[ab\bar\mu]+(|a+b|^2-|a|^2-|b|^2)|\mu|^2\big\} \\
&\quad -\big(2s\Im[a\bar\mu]-2\Re[a\bar\mu]+|a|^2\big)
\big(2s\Im[b\bar\mu]-2\Re[b\bar\mu]+|b|^2\big) \\
%&=\nu^2(1+s^2)\Re[abe^{-i\lambda}/\nu+2a\bar b] \\
%&\quad -\nu^2\big(2s\Im[ae^{-i\lambda}]-2\Re[ae^{-i\lambda}]+|a|^2/\nu\big)
%\big(2s\Im[be^{-i\lambda}]-2\Re[be^{-i\lambda}]+|b|^2/\nu\big) \\
&=ls^2-2ms+n,
\end{align*}
where
\begin{align*}
l&=\Re[ab\bar\mu]+2\Re\big[a\bar\mu\overline{b\bar\mu}\big]-4\Im[a\bar\mu]\Im[b\bar\mu] \\
&=\Re[ab\bar\mu]+2\Re[ab\bar\mu^2]=\Re[ab\bar\mu(1+2\bar\mu)]=\nu L, \\
m&=\Im[a\bar\mu](|b|^2-2\Re[b\bar\mu])+\Im[b\bar\mu](|a|^2-2\Re[a\bar\mu]) \\
&=\Im[ab\bar\mu(\bar a+\bar b-2\bar\mu)]=\nu M, \aand \\
n&=\Re[ab\bar\mu]+2\Re\big[a\bar\mu\overline{b\bar\mu}\big]
-(|a|^2-2\Re[a\bar\mu])(|b|^2-2\Re[b\bar\mu]) \\
&=\Re\big[ab\bar\mu-2ab\bar\mu^2+2ab\bar b\bar\mu+2a\bar ab\bar\mu\big]-|a|^2|b|^2
=\nu N.
\end{align*}
Now the assertion follows from Theorem \ref{thm:main}.
\end{pf}

When $e^{-i\lambda}ab$ or $ab$ is real, the conditions in the theorem are simplified as follows.

\begin{cor}\label{cor:1}
Let $\alpha$ and $\lambda$ be real numbers with $0\le\alpha<1$ and $0<|\lambda|<\pi/2$.
Let $a,b$ be complex numbers with $a+b\ne -1,-2,-3,\dots$
and suppose that $m=e^{-i\lambda}ab$ is a positive real number.
Then the shifted hypergeometric function $z\Gauss(a,b;a+b+1;z)$
is $\lambda$-spirallike of order $\alpha$ if
\begin{align*}
&(\Im[a+b]-(1-\alpha)\sin2\lambda)^2\\
\le~ &\big(2-\alpha+(1-\alpha)\cos2\lambda\big)\big\{2\Re[a+b]+\alpha-(1-\alpha)\cos2\lambda-m/((1-\alpha)\cos\lambda)\big\}.
\end{align*}
\end{cor}

\begin{pf}
In order to apply Theorem \ref{thm:spl},
under the assumptions of the corollary, we compute
\begin{align*}
L&=m[2-\alpha+(1-\alpha)\cos2\lambda], \\
M&=-m(\Im[a+b]-(1-\alpha)\sin2\lambda),
\aand \\
N&=m[2\Re[a+b]+\alpha-(1-\alpha)\cos2\lambda-m/((1-\alpha)\cos\lambda)].
\end{align*}
Thus the assertion follows.
\end{pf}

We note that $L>0$ and $M^2\le LN$ both implies that $N\ge0$ in the assumption
of Theorem \ref{thm:spl}.
We obtain the next result by keeping it in mind.

\begin{cor}\label{cor:2}
Let $\alpha$ and $\lambda$ be real numbers with $0\le\alpha<1$ and
$$
\frac{1-2\alpha}{4(1-\alpha)}< \cos^2\lambda<1.
$$
Let $a,b$ be complex numbers with $a+b\ne -1,-2,-3,\dots$
and suppose that $m=ab$ is a positive real number.
Then the shifted hypergeometric function $z\Gauss(a,b;a+b+1;z)$
is $\lambda$-spirallike of order $\alpha$ if
\begin{align*}
&\left(\frac{\Im[e^{i\lambda}(a+b)]}{\cos\lambda}-2(1-\alpha)\sin2\lambda\right)^2\le
\left(4(1-\alpha)\cos^2\lambda+2\alpha-1\right) \\
&\times\left(\frac{2\Re[e^{i\lambda}(a+b)]}{\cos\lambda}-
4(1-\alpha)\cos^2\lambda+(3-2\alpha)-\frac{m}{(1-\alpha)\cos^2\lambda}\right)
\notag
\end{align*}
\end{cor}

\begin{pf}
Similarly, assuming that $m=ab$ is positive, we compute
\begin{align*}
L&=m\Re\big[e^{-i\lambda}(2-\alpha+(1-\alpha)e^{-2i\lambda})\big]
=m\cos\lambda[4(1-\alpha)\cos^2\lambda-(1-2\alpha)], \\
M&=m\Im\big[e^{-i\lambda}(\bar a+\bar b-2(1-\alpha)e^{-i\lambda}\cos\lambda)\big]\\
&=-m[\Im[e^{i\lambda}(a+b)]-2(1-\alpha)\sin2\lambda\cos\lambda], \aand \\
N&=m\Re\big[e^{-i\lambda}(2\bar a+2\bar b+\alpha-(1-\alpha)e^{-2i\lambda})\big]
-m^2/((1-\alpha)\cos\lambda) \\
&=m\big\{2\Re[e^{i\lambda}(a+b)]+\alpha\cos\lambda-(1-\alpha)\cos3\lambda
-m/((1-\alpha)\cos\lambda)\big\}.
\end{align*}
Here, we used the formula $\cos3\lambda=\cos\lambda(4\cos^2\lambda-3).$
Observe that $L>0$ precisely if $\cos^2\lambda>(1-2\alpha)/4(1-\alpha).$
The assertion now follows from Theorem \ref{thm:spl}.
\end{pf}

When $\alpha=0$, Theorem \ref{thm:spl} and Corollaries \ref{cor:1} and \ref{cor:2} reduce to Theorem 1.4 and Corollaries 1.5 and 1.6 in \cite{SW16sp}, respectively.

\section{Strong starlikeness}

The authors gave in the previous paper \cite{SW16sp} a sufficient condition
for the shifted hypergeometric function to be strongly starlike as in the following.

{\bf Theorem A} (\cite{SW16sp}).
{\itshape
Let $1/3<\alpha<1$ and $a,b$ be complex numbers with
$a+b\in\R$ and $ab>0.$
Then the shifted
hypergeometric function $z\Gauss(a,b;a+b+1;z)$ is strongly starlike
of order $\alpha$ if
$$
\big\{(a-b)^2+6(a+b)-3\big\}\sin^2\frac{\pi\alpha}2\ge a^2+ab+b^2.
$$
}

However, this was obtained as a corollary of the spirallikeness result
(Theorem \ref{thm:spl} with $\alpha=0$).
Therefore, the unpleasant assumption $1/3<\alpha<1$ was inevitable.
The following can be obtained as a consequence of the main theorem.
Though the condition is much involved, this restriction does not appear explicitly.

\begin{thm}\label{thm:ss}
Let $\alpha$ be a real number with $0<\alpha<1$ and let $a,b,c$ be
complex numbers satisfying $ab\ne0$ and $c\ne0,-1,-2,\dots$.
The shifted hypergeometric function $f(z)=z\Gauss(a, b; c;z)$ is strongly starlike
of order $\alpha$ if the following three conditions are satisfied:
\renewcommand{\labelenumi}{({\roman{enumi}})}
\begin{enumerate}
\item
$p=a+b+1-c$ is a real number,
\item
$|\arg(ab-p)|<\pi\alpha/2,$
\item
$
G_\varepsilon(s^\alpha)\leq
\alpha(s+s\inv)s^{\alpha}\Re\left[(ab-p)e^{\varepsilon\pi i(1-\alpha)/2}\right],$
for $s\in(0,+\infty),~\varepsilon=\pm 1,$
where $G_\varepsilon(x)=Sx^3+T_\varepsilon x^2+U_\varepsilon x+V$ with
\begin{align*}
S&=2p\cos\frac{\pi\alpha}2, \\
T_\varepsilon&=|a|^2+|b|^2-|c-1|^2-2p-4p\cos^2\frac{\pi\alpha}2
+4\Re[ae^{-\varepsilon i\pi\alpha/2}]\Re[be^{-\varepsilon i\pi\alpha/2}], \\
U_\varepsilon&=-2(|a|^2+|b|^2-|c-1|^2-3p)\cos\frac{\pi\alpha}2 \\
&\quad +2\Re[ae^{-\varepsilon i\pi\alpha/2}](|b|^2-2\Re b)
+2\Re[be^{-\varepsilon i\pi\alpha/2}](|a|^2-2\Re a), \\
V&=|a|^2+|b|^2-|c-1|^2-2p+(|a|^2-2\Re a)(|b|^2-2\Re b).
\end{align*}
\end{enumerate}
\end{thm}

\begin{pf}
To prove that $f\in\es^*(\phi_\alpha),$
it is sufficient to check the conditions in Theorem \ref{thm:main} for
$$
Q(z)=\phi_{\alpha}(z)-1=\left(\frac{1+z}{1-z}\right)^{\alpha}-1
$$
with $\alpha\in[0,1)$.

Let $\zeta=e^{i\theta}$ for $0<|\theta|<\pi,$ $\varepsilon=\theta/|\theta|=\sgn\theta,$
and $s=\cot(|\theta|/2)\in (0,+\infty).$
As in the previous section, we have
\begin{align*}
%\zeta&=e^{i\theta}=\frac{s^2-1}{s^2+1}+i\frac{2s}{s^2+1},\\
Q(\zeta)&
%=\left(\frac{1+e^{i\theta}}{1-e^{i\theta}}\right)^{\alpha}-1
=(\varepsilon is)^{\alpha}-1
=e^{i\beta}s^\alpha-1
\aand \\
\zeta Q'(\zeta)&=[Q(\zeta)+1]\frac{2\alpha e^{i\theta}}{1-e^{2i\theta}}
=-\frac{\alpha e^{i\beta}s^\alpha(1+s^2)}{2i\varepsilon s}
=-\frac{\alpha}2 e^{-i\gamma}s^\alpha(s+s\inv),
\end{align*}
where
$$
\beta=\varepsilon\frac{\pi\alpha}2
\aand
\gamma=\varepsilon\frac{\pi(1-\alpha)}2.
$$

Hence by using these relations, we get
\begin{align*}
-2\Re[(pQ(\zeta)+ab)\overline{\zeta Q'(\zeta)}]
&=\Re\left[(pe^{i\beta}s^\alpha+ab-p)\alpha e^{i\gamma}s^\alpha(s+s\inv)\right]\\
&=\alpha(s+s\inv)s^\alpha\big\{-\varepsilon s^{\alpha}\Im p
+\Re\left[e^{i\gamma}(ab-p)\right]\big\}.\\
\end{align*}

By \eqref{eq:main1}, we need the conditions
$$
\Im p=0\quad \text{and}\quad \Re\left[e^{i\gamma}(ab-p)\right]>0
$$
for $\varepsilon=\pm1.$
In this way, we arrived at the first and second conditions in the theorem.
From now on, we assume that these two conditions are satisfied.
In view of \eqref{eq:AB}, we compute
\begin{align*}
&\quad ~ |B(\zeta)|^2-|A(\zeta)|^2\\
&=(s^{2\alpha}-2s^\alpha\cos\beta+1)(2ps^\alpha\cos\beta-2p+|a|^2+|b|^2-|c-1|^2) \\
&\quad +(2\Re[ae^{-i\beta}s^\alpha-a]+|a|^2)(2\Re[be^{-i\beta}s^\alpha-b]+|b|^2) \\
&=G_\varepsilon(s^\alpha).
\end{align*}
Therefore the condition \eqref{eq:main2} can be
presented as in the third condition of the theorem.
\end{pf}

We now assume that $p=a+b+1-c=0.$
Letting $\eta=e^{-\varepsilon i\pi\alpha/2},$ we compute
the coefficients of $G_\varepsilon(x)$ in Theorem \ref{thm:ss}
by making use of \eqref{eq:zw} as follows: $S=0,$
\begin{align*}
T_\varepsilon&=-2\Re[a\bar b]+4\Re[a\eta]\Re[b\eta]=2\Re[ab\eta^2], \\
U_\varepsilon&=2\Re[a\bar b](\eta+\bar\eta)+2\Re[a\eta(b\bar b-b-\bar b)
+b\eta(a\bar a-a-\bar a)] \\
&=2\Re[ab\eta(\bar a+\bar b-2)], \\
V&=-2\Re[a\bar b]+|ab|^2-2\Re[ab\bar b+a\bar ab]+4\Re a\Re b \\
&=2\Re[ab]+|ab|^2-2\Re[ab(\bar a+\bar b)].
\end{align*}
Thus, we obtain the following corollary.

\begin{cor}\label{cor:sst}
Let $\alpha$ be a real number with $0<\alpha<1$ and let $a,b$ be complex numbers
with $a+b\ne -1,-2,-3,\dots.$
The shifted hypergeometric function $z\Gauss(a, b; a+b+1;z)$ is strongly starlike
of order $\alpha$ if
\renewcommand{\labelenumi}{({\roman{enumi}})}
\begin{enumerate}
\item
$|\arg(ab)|<\pi\alpha/2,$ and
\item
\begin{align*}
&\qquad 2s^{2\alpha}\Re[e^{-\varepsilon\pi i\alpha}ab]
+2s^\alpha \Re[e^{-\varepsilon\pi i\alpha/2}ab(\bar a+\bar b-2)] \\
&\quad +|ab|^2-2\Re[ab(\bar a+\bar b-1)]
\leq \alpha(s+s\inv)s^{\alpha}\Re\big[e^{\varepsilon\pi i(1-\alpha)/2}ab\big]
\end{align*}
for all $s\in(0,+\infty)$ and both signs $\varepsilon=\pm1.$
\end{enumerate}
\end{cor}

Since the condition in the last corollary is still implicit,
we make a crude estimate.
We first prepare the following elementary lemma.

\begin{lem}
Let $\alpha, A,B,C$ and $K$ be constants with $0<\alpha<1$ and $K>0.$
If $B/2+\max\{A,C\}\le K,$ then
$$
As^\alpha+B+Cs^{-\alpha}\le K(s+s\inv)
\quad \text{for}~ s\in(0,+\infty).
$$
\end{lem}

\begin{pf}
It is easy to observe that $s^\alpha+s^{-\alpha}$ is increasing in
$0\le \alpha\le 1$ for a fixed $s>0, s\ne1.$
In particular, $2<s^\alpha+s^{-\alpha}<s+s\inv$ for $0<\alpha<1, s>0, s\ne1.$
Thus the assertion follows.
\end{pf}

We can now deduce the following from Corollary \ref{cor:sst}.

\begin{cor}
Let $\alpha$ be a real number with $0<\alpha<1.$
Assume that complex numbers $a,b$ with $a+b\ne -1,-2,-3,\dots$ satisfy
the conditions:
\renewcommand{\labelenumi}{({\roman{enumi}})}
\begin{enumerate}
\item
$|\arg(ab)|<\pi\alpha/2,$ and
\item
$\displaystyle
\Re[e^{-\varepsilon\pi i\alpha/2}ab(\bar a+\bar b-2)]
+\max\big\{2\Re[e^{-\varepsilon\pi i\alpha}ab],
|ab|^2-2\Re[ab(\bar a+\bar b-1)]\big\}
$

$\displaystyle
\le \alpha\Re\big[e^{\varepsilon\pi i(1-\alpha)/2}ab\big],
\quad \text{for}~\varepsilon=\pm1.
$
\end{enumerate}
Then the function $z\Gauss(a, b; a+b+1;z)$ is strongly starlike
of order $\alpha.$
\end{cor}

We further assume that $l=a+b$ is real and $m=ab$ is real and positive.
Then the condition (ii) in the last corollary reads
$$
(l-2)\cos\frac{\pi\alpha}2+\max\{2\cos\pi\alpha, m-2(l-1)\}
\le \alpha\sin\frac{\pi\alpha}2.
$$
In particular, if $2\cos\pi\alpha\ge m-2(l-1),$
or equivalently, if $m-2l+4=(a-2)(b-2)\le 2\cos\pi\alpha+2
=4\cos^2(\pi\alpha/2),$ the condition takes the form $(l-2)\cos(\pi\alpha/2)
+2\cos\pi\alpha\le\alpha\sin(\pi\alpha/2).$
Therefore, we finally obtain the following corollary.

\begin{cor}
Let $\alpha$ be a real number with $0<\alpha<1$ and assume that
$a,b$ are complex numbers such that $a+b$ is real and $ab>0.$
If, in addition,
$$
(a-2)(b-2)\le 4\cos^2\frac{\pi\alpha}2
\aand
a+b\le 2-2\frac{\cos\pi\alpha}{\cos\frac{\pi\alpha}{2}}+
\alpha\tan\frac{\pi\alpha}2,
$$
then the function $z\Gauss(a, b; a+b+1;z)$
is strongly starlike of order $\alpha.$
\end{cor}

We remark that, under the hypothesis in the last corollary, $l=a+b$ should satisfy
the inequalities
$$
2\sin^2\frac{\pi\alpha}2<\frac{ab}2+2\sin^2\frac{\pi\alpha}2\le l
\le 2-2\frac{\cos\pi\alpha}{\cos\frac{\pi\alpha}{2}}+\alpha\tan\frac{\pi\alpha}2.
$$
We can see that $\big[2-2\cos\pi\alpha/\cos(\pi\alpha/2)+\alpha\tan(\pi\alpha/2)\big]
-2\sin^2(\pi\alpha/2)$ is positive for $0<\alpha<1.$
Therefore, there are some $a,b\in\C$ satisfying the hypothesis in the corollary
for each $0<\alpha<1.$
Compare with Theorem A.

\medskip

\noindent
{\bf Acknowledgements.}
The authors sincerely thank the referees for corrections and helpful suggestions.

\def\cprime{$'$} \def\cprime{$'$} \def\cprime{$'$}
\providecommand{\bysame}{\leavevmode\hbox to3em{\hrulefill}\thinspace}
\providecommand{\MR}{\relax\ifhmode\unskip\space\fi MR }
\providecommand{\MRhref}[2]{%
  \href{http://www.ams.org/mathscinet-getitem?mr=#1}{#2}
}
\providecommand{\href}[2]{#2}

%\bibliography{papers}

\begin{thebibliography}{10}

\bibitem{AS:hand}
M.~Abramowitz and I.~A. Stegun, \emph{Handbook of {M}athematical {F}unctions},
  Dover, 1972.

\bibitem{Ahlfors:ca}
L.~V. Ahlfors, \emph{Complex {A}nalysis, 3rd ed.}, McGraw Hill, New York, 1979.

\bibitem{BPV02}
 R.~Balasubramanian, S.~Ponnusamy and M.~Vuorinen, \emph{On hypergeometric functions and function spaces}, J. Comput. Appl. Math. (2) \textbf{139} (2002), 299--322.

\bibitem{BNP06}
R.~W.~Barnard, S.~Naik and S.~Ponnusamy, \emph{Univalency of weighted integral transforms of certain functions}, J. Comput. Appl. Math. \textbf{193} (2006), 638--651.

\bibitem{BK69}
D.~A. Brannan and W.~E. Kirwan, \emph{On some classes of bounded univalent
  functions}, J. London Math. Soc. (2) \textbf{1} (1969), 431--443.



\bibitem{FKZ76}
M.~Fait, J.~G. Krzy\.z, and J.~Zygmunt, \emph{Explicit quasiconformal
  extensions for some classes of univalent functions}, Comment. Math. Helv.
  \textbf{51} (1976), 279--285.

\bibitem{Goodman:univ}
A.~W. Goodman, \emph{Univalent {F}unctions, 2 vols.}, Mariner Publishing Co.
  Inc., 1983.

\bibitem{HPV10}
P.~H\"asto, S.~Ponnusamy, and M.~Vuorinen, \emph{Starlikeness of the {G}aussian
  hypergeometric functions}, Complex Var. Elliptic Equ. \textbf{55} (2010),
  173--184.

\bibitem{Kus02}
R.~K\"ustner, \emph{Mapping properties of hypergeometric functions and
  convolutions of starlike or convex functions of order $\alpha$}, Comput.
  Methods Funct. Theory \textbf{2} (2002), 597--610.

\bibitem{Kus07}
R.~K{\"u}stner, \emph{On the order of starlikeness of the shifted {G}auss
  hypergeometric function}, J. Math. Anal. Appl. \textbf{334} (2007),
  1363--1385.

\bibitem{Lewis79}
J.~Lewis, \emph{Applications of a convolution theorem to {J}acobi polynomials},
  SIAM J. Math. Anal. \textbf{10} (1979), 1110--1120.

\bibitem{Libera67}
R.~J. Libera, \emph{Univalent $\alpha$-spiral functions}, Can. J. Math.
  \textbf{19} (1967), 449--456.

\bibitem{MM92B}
W.~Ma and D.~Minda, \emph{A unified treatment of some special classes of
  univalent functions}, Proceedings of the {C}onference on {C}omplex {A}nalysis
  (Z.~Li, F.~Ren, L.~Yang, and S.~Zhang, eds.), International Press Inc., 1992,
  pp.~157--169.

\bibitem{MM:ds}
S.~S. Miller and P.~T. Mocanu, \emph{Differential {S}ubordinations. {T}heory
  and {A}pplications}, Marcel Dekker, Inc., New York, 2000.

\bibitem{Pin68}
B.~Pinchuk, \emph{On starlike and convex functions of order $\alpha$}, Duke
  Math. J. \textbf{35} (1968), 721--734.

\bibitem{Pom:bound}
{Ch}. Pommerenke, \emph{Boundary {B}ehaviour of {C}onformal {M}aps},
  Springer-Verlag, 1992.

\bibitem{Rob36}
M.~S. Robertson, \emph{On the theory of univalent functions}, Ann. of Math.
  \textbf{37} (1936), 374--408.

\bibitem{Rus:conv}
{St}. Ruscheweyh, \emph{Convolutions in {G}eometric {F}unction {T}heory},
  S\'eminaire de Math\'ematiques Sup\'erieures, vol.~83, Les Presses de
  l'Universit\'e de Montr\'eal, Montr\'eal, 1982.

\bibitem{Stank66}
J.~Stankiewicz, \emph{Quelques probl\`emes extr\'emaux dans les classes des
  fonctions $\alpha$-angulairement \'etoil\'ees}, Ann. Univ. Mariae Curie-Sk\l
  odowska, Sect. A \textbf{20} (1966), 59--75.

\bibitem{SugawaDual}
T.~Sugawa, \emph{A self-duality of strong starlikeness}, Kodai Math J.
  \textbf{28} (2005), 382--389.

\bibitem{SW16sp}
T.~Sugawa and L.-M. Wang, \emph{Spirallikeness of shifted hypergeometric
  functions}, to appear in Ann. Acad. Sci. Fenn. Math.

\bibitem{Temme:sp}
N.~M. Temme, \emph{Special {F}unctions}, John Wiley \& Sons, Inc., New York,
  1996.

\bibitem{WW:anal}
E.~T. Whittaker and G.~N. Watson, \emph{A course of {M}odern {A}nalysis. {A}n
  {I}ntroduction to the {G}eneral {T}heory of {I}nfinite {P}rocesses and of
  {A}nalytic {F}unctions; with an {A}ccount of the {P}rincipal {T}ranscendental
  {F}unctions}, Fourth edition. Reprinted, Cambridge University Press, New
  York, 1962.

\end{thebibliography}
\end{document}